\documentclass[11pt,english]{smfart}
\usepackage{graphpap}
\usepackage{amssymb}
\newtheorem{theoreme}{Theorem}
\renewcommand{\theequation}{\arabic{equation}}
\newtheorem{proposition}{Proposition}
\newtheorem{lemme}[proposition]{Lemma}

\newtheorem{definition}[proposition]{Definition}
\newtheorem{remarque}[proposition]{Remark}

\numberwithin{equation}{section}
\numberwithin{proposition}{section}

\def\Re{\textrm{Re}} 
\def\11{{\rm 1~\hspace{-1.4ex}l} }
\def\R{\mathbb R}
\def\C{\mathbb C}

\def\N{\mathbb N}

\def\O{\mathbb O}

\begin{document}
\title
[Supercritical wave equations]
{Random data Cauchy theory for supercritical wave equations II : a global existence result}
\author{Nicolas Burq}
\address{D\'epartement de Math\'ematiques, Universit\'e Paris XI, 91 405  Orsay Cedex, France, 
and Institut Universitaire de France}
\email{nicolas.burq@math.u-psud.fr}
\author{Nikolay Tzvetkov}
\address{D\'epartement de Math\'ematiques, Universit\'e Lille I, 59 655 Villeneuve d'Ascq Cedex, France}
\email{nikolay.tzvetkov@math.univ-lille1.fr}
\begin{abstract} 
We prove that the subquartic wave equation on the three dimensional ball $\Theta$, with Dirichlet boundary conditions 
admits global strong solutions for a large set of random supercritical initial
data in $\cap_{s<1/2} H^s( \Theta)$. 
We obtain this result 
as a consequence of a general random data Cauchy theory for supercritical wave
equations developed in our previous work~\cite{BT2} 
and invariant measure considerations which allow us to obtain also precise
large time dynamical informations on our solutions. 
\end{abstract}
\maketitle
%
\section{Introduction}
In our previous work~\cite{BT2}, we developed a general theory for constructing local strong solutions to nonlinear 
wave equations, 
posed on compact riemannian manifolds with supercritical random initial data.
The goal of this article is to show that in a {\it very particular case} we
can combine this local theory with some 
invariant measure arguments (see the work by Bourgain~\cite{Bo1, Bo2} and the authors~\cite{Tz1, Tz2, BT1}) to obtain {\em global} solutions. Namely, we shall
consider the nonlinear wave equation 
with Dirichlet boundary condition  posed on $\Theta$, the unit ball of $\R^3$,
\begin{equation}\label{1}
(\partial_{t}^{2}-\mathbf{\Delta})w+|w|^{\alpha}w=0,\quad (w,\partial_t w)|_{t=0}=(f_1,f_2), 
\quad u\mid_{\R_t \times \partial \Theta} =0,\qquad \alpha>0
\end{equation}
with radial {\it real valued} initial data $(f_1,f_2)$. Our aim is to give a proof of the following result.
\begin{theoreme}\label{thm1}
Suppose that $\alpha<3$. Let us fix a real number $p$ such that $\max(4,2\alpha)<p<6$.
Let $((h_n(\omega),l_n(\omega))_{n=1}^{\infty}$ be a sequence of independent
standard real Gaussian random variables on a probability space
$(\Omega,{\mathcal A},p)$. Consider (\ref{1}) 
with initial data 
$$
f_1^\omega(r)=\sum_{n=1}^{\infty}\frac{ h_{n}(\omega)}{n\pi}e_{n}(r),\quad
f_2^\omega(r)=\sum_{n=1}^{\infty} l_{n}(\omega)e_{n}(r)\,,
$$
where $(e_n(r))_{n=1}^{\infty}$ is the orthonormal basis consisting in  
radial eigenfunctions of the Laplace operator
with Dirichlet boundary conditions, associated to eigenvalues $(\pi n)^2.$

Then for every $s<1/2$, almost surely in $\omega\in \Omega$, the problem (\ref{1}) has a unique global solution 
$$
u^\omega\in C(\R,H^s(\Theta))\cap L^p_{loc}( \mathbb{R}_t;L^p(\Theta))\,.
$$
Furthermore, the solution is a perturbation of the {\em linear } solution
$$ 
u^\omega (t) = U(t)(f_1^\omega, f_2^\omega)+ v^\omega(t)= 
\cos (\sqrt{- \mathbf{\Delta}}t)f_1^\omega + \frac{ \sin{\sqrt {\mathbf{-\Delta}}}}{\sqrt{ - \mathbf{\Delta}}}
f_2^\omega+ v^\omega (t),
$$ 
where $v^\omega\in C(\R,H^\sigma(\Theta))$ for some $\sigma>1/2$.
Moreover
$$
\|u^\omega(t)\|_{H^s(\Theta)} \leq 
C(\omega,s) \log(2+|t|)^{\frac 1 2}\,.
$$
\end{theoreme}

\begin{remarque}{\rm
We have  (see~\cite[Lemma 3.2]{BT1} that almost surely 
$$
(f_1^\omega,f_2^\omega)\in \bigcap_{s<1/2}(H^{s}(\Theta)\times H^{s-1}(\Theta))
$$
but the probability of the event 
$\{(f_1^\omega,f_2^\omega)\in H^{\frac{1}{2}}(\Theta)\times H^{-\frac{1}{2}}(\Theta)\}$ 
is zero. Thus in the above statement, we obtain global solutions for data which are not in 
$H^{1/2}(\Theta)\times H^{-1/2} (\Theta)$. Observe that the equation~\eqref{1} is $H^{3/2 - 2/\alpha}\times H^{1/2- 2/\alpha}$ critical. 
As a consequence, for $2< \alpha <3$, we obtain global existence of strong solutions for a supercritical model, 
a result which seems to be completely out of reach of the present deterministic methods, 
even for the local existence theory.} 
\end{remarque}
\begin{remarque} {\rm 
The initial data we consider can be arbitrarily large in $L^2( \Theta)$. Consequently, our result is {\em not} a ``small data  result''.}
\end{remarque}

The rest of the paper is organized as follows. In the next section, we analyze
the Hamiltonian structure of the equation and we introduce the suitable finite
dimensional model which approximate (\ref{1}). In Section~3, we establish the
key probabilistic estimate concerning the $L^p$ space-time norms of the free evolution.
In Section~4, we recall the deterministic Strichartz estimates for the free
evolution. Section~5 is devoted to local well-posedness results. 
In Section~6, we perform the globalization argument. The analysis of Section~6 has much in
common with some arguments already appeared in \cite{Tz1,Tz2,BT1}.
\section{Reduction of the problem and approximating ODE}
For $\sigma\in\R$, we define (see~\cite{BT1} for more details) $H^{\sigma}(\Theta)$ as
$$
H^{\sigma}(\Theta)\equiv\Big\{u=\sum_{n=1}^{\infty}c_{n}e_{n},\,c_{n}\in\C\,:\,
\|u\|_{H^\sigma( \Theta)}^2= \sum_{n=1}^{\infty} (n\pi)^{2\sigma}|c_n|^2<\infty\Big\}
$$
Remark (see ~\cite{LionsMagenes}) that for $-1/2<\sigma<1/2$ these spaces coincide with the classical Sobolev 
spaces (of radial functions) and are independent of the choice of boundary conditions we made.
Following~\cite{BT1}, we make some algebraic manipulations on (\ref{1}) allowing 
to write it as a first order
equation in $t$. Set $u\equiv w+i{\sqrt{ - \mathbf{\Delta}}^{-1}}{\partial_t w}$. 
Then we have that $u$ solves the equation  
\begin{equation}\label{2}
(i\partial_{t}-\sqrt{ - \mathbf{\Delta}})u-\sqrt{ - \mathbf{\Delta}}^{-1}\big(|\Re(u)|^{\alpha}\Re(u)\big)
=0,\quad u|_{t=0}=u_{0},\quad u|_{\R\times\partial\Theta}=0,
\end{equation}
where $u_0=f_1+i\sqrt{ - \mathbf{\Delta}}^{-1}f_2$.
Equation~\eqref{2} is (formally) an Hamiltonian equation on $L^2(\Theta)$ with Hamiltonian, 
$$
\frac{1}{2}\|\sqrt{ - \mathbf{\Delta}} (u)\|_{L^2(\Theta)}^{2}+
\frac{1}{\alpha+2}\|\Re(u)\|_{L^{\alpha+2}(\Theta)}^{\alpha+2}
$$
which is (formally) conserved by the flow of (\ref{2}).

In order to prove Theorem~\ref{thm1}, we will need to study (\ref{2}) with initial data given by
$$
u_{0}(r,\omega)=\sum_{n=1}^{\infty}\frac{ g_{n}(\omega)}{n\pi}e_{n}(r),
$$
where $g_{n}(\omega)=h_{n}(\omega)+il_{n}(\omega)$ are independent normalized complex Gaussian.

For $N\geq 1$, we denote by $E_{N}$ the $N$ dimensional vector space on $\C$ spanned by $(e_n)_{n=1}^{N}$.
Fix $\chi \in C^\infty_0(-1,1)$ equal to $1$ on $(-1/2, 1/2)$. 
Let us define $S_{N}=\chi(\frac {- \mathbf{\Delta} }{\pi^2 N^2})$. 
This  operator sends $L^2$ to  $E_{N}$  and satisfies
$$
S_{N}\Big(\sum_{n=1}^{\infty}c_{n}e_{n}\Big)=\sum_{n=1}^{\infty}\chi(\frac {n^2} { N^2})c_{n}e_{n}\,.
$$
Let us observe that the map $S_N$ we use in this paper is slightly different than the one involved 
in \cite{BT1}. The reason is that we will use $L^p$, $p\neq 2$ mapping properties of $S_N$ which do not 
hold for the map used in \cite{BT1}.
More precisely, we have the following statement.
\begin{lemme}\label{psdo}
Let us fix $p\in [1,\infty]$. There exists $C>0$ such that for every integer
$N\geq 1$, $\|S_{N}\|_{L^p\rightarrow L^p}\leq C$. Moreover, for every $f\in
L^p$, $S_{N}(f)\rightarrow f$ as $N\rightarrow\infty$ in $L^p$.
\end{lemme}
The proof of Lemma~\ref{psdo} is essentially contained in \cite{BGT}, where
the case of boundary less manifold is considered. The general case is more
involved technically and requires a precise description of the operator $S_N$
(see for example~\cite[section 4.3]{Sz}). 
Notice however that the radial assumption we made here would allow a rather direct proof.

We shall approximate the solutions of (\ref{2}) by the solutions of the ODE
\begin{equation}\label{pak_N}
(i\partial_{t}-\sqrt{ - \mathbf{\Delta}})u-
S_{N}\left(\sqrt{ - \mathbf{\Delta}}^{-1}\big(|S_N(\Re(u))|^{\alpha}S_N(\Re(u))\big)\right)=0,\,
u|_{t=0}=u_{0}\in E_{N}\,\,.
\end{equation}
Let us define the measure $\mu_{N}$ on $E_{N}$ as the image measure under the map from
$(\Omega,{\mathcal A},p)$ to $E_{N}$ (equipped with the Borel sigma algebra) defined by
\begin{equation}\label{zve}
\omega\longmapsto \sum_{n=1}^{N}\frac{h_{n}(\omega)+i l_{n}(\omega)}{n\pi}e_{n}\,,
\end{equation}
where $h_{n}(\omega),l_{n}(\omega)$, $n=1,\cdots N$  is a sequence of
independent standard real gaussians ($h_n, l_n\in {\mathcal N}(0,1)$).
We next define the measure $\rho_{N}$ as the image measure by the map
(\ref{zve}) of the measure 
\begin{equation*} 
\exp
\Big(-\frac{1}{\alpha+2}\|S_N(\sum_{n=1}^{N}\frac{ h_{n}(\omega)}{n\pi}e_{n})\|_{L^{\alpha+2}(\Theta)}^{\alpha+2}\Big)
dp(\omega)\, .
\end{equation*}
It turns out that $\rho_N$ is invariant under the flow of (\ref{pak_N}).
\begin{proposition}\label{liouville}
For every $u_0\in E_N$ the flow of (\ref{pak_N}) is defined globally in time.
Moreover the measure $\rho_N$ is invariant under this flow.
\end{proposition}
\begin{proof}
The proof of this result is essentially (in a slightly different context) in~\cite{BT1}. 
For the sake of completeness, we recall briefly the proof. 
The local existence and uniqueness for the ODE (\ref{pak_N}) follows from the
Cauchy-Lipschitz theorem. We can 
extend globally in time the solutions of (\ref{pak_N}) thanks to the energy
conservation law associated to (\ref{pak_N}). Indeed if we multiply (\ref{pak_N})
by $\mathbf{\Delta} u-S_{N}(S_N(|\Re(u))|^{\alpha}S_N(\Re(u)))$ (which is an element of $E_N$) 
and integrate over $\Theta$, we get that the solutions of (\ref{pak_N}) satisfy
$$
\frac{d}{dt}
\Big[\frac{1}{2}\|\sqrt{ - \mathbf{\Delta}} (u)\|_{L^2(\Theta)}^{2}+
\frac{1}{\alpha+2}\|S_N(\Re(u))\|_{L^{\alpha+2}(\Theta)}^{\alpha+2}\Big]=0\,.
$$
Thus we have a control uniform with respect to time  and therefore the solutions of (\ref{pak_N}) 
are defined globally in time. Let us now turn to the proof of the measure invariance.
Let us decompose the solution of (\ref{pak_N}) as
$$
u(t)=\sum_{n=1}^{N}\big(a_{n}(t)+ib_{n}(t)\big)e_{n},\quad a_{n}(t),b_n(t)\in\R\,.
$$
Then, if we set
$$
H(a_1,\dots,a_N,b_1,\dots,b_N)\equiv
\frac{1}{2}\sum_{n=1}^{N}(\pi n)^2(a_n^2+b_n^2)+\frac{1}{\alpha+2}\int_{\Theta}
\big|S_N(\sum_{n=1}^{N}a_{n}e_{n})\big|^{\alpha+2}
$$
the problem (\ref{pak_N}) may be rewritten in the coordinates $a_n,b_n$ as
\begin{equation}\label{ab}
\dot{a}_{n}=(\pi n)^{-1}\frac{\partial H}{\partial b_n},\quad 
\dot{b}_{n}=-(\pi n)^{-1}\frac{\partial H}{\partial a_n},\quad n=1,\dots,N
\end{equation}
($e_n$ are real valued). Let us first observe that thanks to the structure of
(\ref{ab}) the quantity $H(a_1,\dots,a_N,b_1,\dots,b_N)$
is conserved under the flow of (\ref{ab}). Therefore we may apply Liouville's theorem for 
divergence free vector fields to obtain that the measure 
$$
\prod_{n=1}^{N}(\pi n)^2da_n db_n
$$
is conserved by the flow of (\ref{ab}).
Since $H(a_1,\dots,a_N,b_1,\dots,b_N)$ is also conserved under the flow of
(\ref{ab}) we obtain that the measure
\begin{multline*}
\frac 1 {(2\pi)^N}\exp\big(-H(a_1,\dots,a_N,b_1,\dots,b_N)\big)\prod_{n=1}^{N}(\pi n)^2 da_n
db_n
\\
=
\exp\Big(-\frac{1}{(\alpha+2)}\int_{\Theta}
\big|S_N(\sum_{n=1}^{N}a_{n}e_{n})\big|^{\alpha+2}\Big)
\prod_{n=1}^{N}
\sqrt{\frac{\pi} 2}  n  e^{-(\pi n)^2(a_n^2/2)}da_n \sqrt{\frac {\pi} 2} ne^{-(\pi n)^2(b_n^2/2)}db_n
\end{multline*}
is also conserved by the flow of (\ref{ab})
which, coming back to $E_N$, completes the proof of Proposition~\ref{liouville}.
\end{proof}
Let us fix from now on in the rest of this paper a number $s<1/2$.
Let us define the measure $\mu$ on $H^{s}(\Theta)$ as the image measure under the map from
$(\Omega,{\mathcal A},p)$ to $H^{s}(\Theta)$ equipped with the Borel sigma algebra, defined by
\begin{equation}\label{map}
\omega\longmapsto \sum_{n=1}^{\infty}\frac{ h_{n}(\omega)+i l_{n}(\omega)  }{n\pi}e_{n}\,,
\end{equation}
where $((h_{n},l_{n}))_{n=1}^{\infty}$ is a sequence of independent standard real Gaussian random variables. 

Using \cite[Theorem~4]{AT}, we have that for $\alpha<4$ the quantity
$$
\|\sum_{n=1}^{\infty}\frac{ h_{n}(\omega)+i l_{n}(\omega)  }{n\pi}e_{n}\|_{L^{\alpha+2}(\Theta)}
$$
is finite almost surely. 
Therefore, we can define a nontrivial measure $\rho$ as the image measure on
$H^s(\Theta)$ by the map (\ref{map}) of the measure
$$
\exp\Big(-\frac{1}{\alpha+2}
\|\sum_{n=1}^{\infty}\frac{ h_{n}(\omega)}{n\pi}e_{n})\|_{L^{\alpha+2}(\Theta)}^{\alpha+2}\Big)dp(\omega)\,.
$$
Observe that if a Borel set $A\subset H^s(\Theta)$ is of full $\rho$ measure
then $A$ is also of full $\mu$ measure. Therefore, we need to solve (\ref{2}) globally in time
for $u_0$ in a set of full $\rho$ measure.

We next turn to the limits of the measures $\rho_N$. 
We have the following statement.
\begin{lemme}\label{deg}
Set 
\begin{equation}\label{f_N}
f(u)=\exp\Big(-\frac{1}{\alpha+2}\|u\|_{L^{\alpha+2}(\Theta)}^{\alpha+2}\Big),\quad
f_{N}(u)=\exp\Big(-\frac{1}{\alpha+2}\|S_{N}(u)\|_{L^{\alpha+2}(\Theta)}^{\alpha+2}\Big)\, .
\end{equation}
Then
$$
\lim_{N\rightarrow\infty}\int_{H^s(\Theta)}|f_{N}(u)-f(u)|d\mu(u)=0\,.
$$
In particular,
\begin{equation}\label{q}
\lim_{N\rightarrow \infty}\rho_{N}(E_{N})=\rho\big(H^{s}(\Theta)\big).
\end{equation}
\end{lemme}
\begin{proof}
The argument is very close to the proof of \cite[Lemma~3.7]{Tz2} and
therefore we will only sketch it. Thanks to the analysis of \cite{AT} (see also \cite[Lemma~2.3]{Tz2}), we have
that $\|S_{N}(u)\|_{L^{\alpha+2}(\Theta)}^{\alpha+2}$ converges, as $N\rightarrow \infty$ to $\|u\|_{L^{\alpha+2}(\Theta)}^{\alpha+2}$
in $L^1(d\mu)$. Therefore $f_{N}(u)$ converges in measure, as $N\rightarrow \infty$ to
$f(u)$ with respect to the measure $d\mu$. For $\varepsilon >0$, we consider the set
$$
A_{N,\varepsilon}\equiv \big(
u\in H^{s}(\Theta)\,:\, |f_{N}(u)-f(u)|\leq \varepsilon
\big).
$$
Using the Cauchy-Schwarz inequality, we get
$$
\int_{A_{N,\varepsilon}^{c}}|f_{N}(u)-f(u)|d\mu(u)\leq
\|f_N-f\|_{L^{2}(d\mu)}[\mu(A_{N,\varepsilon}^{c})]^{\frac{1}{2}}\leq
2[\mu(A_{N,\varepsilon}^{c})]^{\frac{1}{2}}\,.
$$
On the other hand
$$
\int_{A_{N,\varepsilon}}|f_{N}(u)-f(u)|d\mu(u)\leq\varepsilon\,.
$$
Finally, the convergence in measure of $f_N$ to $f$ implies that for a fixed $\varepsilon$,
$$
\lim_{N\rightarrow \infty}\mu(A_{N,\varepsilon}^{c})=0.
$$
This completes the proof of Lemma~\ref{deg}.
\end{proof}
\section{Gaussian estimates}
Let us recall the following  standard Gaussian estimate (see e.g. \cite{Tz1,Tz2,BT1}).
\begin{lemme}\label{lem.gauss}
Let $c$ be a positive constant satisfying $c<\pi/2$.
Denote by $B(0, \lambda)_s$ the open ball of center $0$ and radius $\lambda$ in $H^s(\Theta)$. 
Then there exists $C_s>0$ 
such that for every $N, \lambda$,
\begin{equation*}
\rho_{N}( B(0,\lambda)_s^c\cap E_N)\leq \mu_N( B(0,\lambda)^c\cap E_N)\leq C_s e^{-c\lambda^2}\,.
\end{equation*} 
\end{lemme}
Let $S(t) = e^{-it\sqrt{ -\mathbf{\Delta}}}$ denote the free evolution operator.
Let us observe that for every $t\in\R$, $S(t+2)=S(t)$.
The following large deviation estimate will play a crucial role in our analysis.
\begin{proposition}\label{largedev}
For any $2\leq p<6$, there exists $C,c>0$ such that for any $N, \lambda>0$,
\begin{equation*}
\mu_N(\{f\in E_N\, :\, \|S(t) f\|_{L^p((0,2)\times \Theta)}>\lambda\})\leq C e^{-c\lambda^2}
\end{equation*} 
\end{proposition}
\begin{proof}
We need to show that there exists $C,c>0$ such that for any $N, \lambda>0$,
\begin{equation}\label{ev1}
p(\{\omega\in \Omega\,:\, \|S(t) f^{\omega}_{N}\|_{L^p((0,2)\times \Theta)}>\lambda\})\leq C e^{-c\lambda^2},
\end{equation}
where
$$
f^{\omega}_{N}=\sum_{n=1}^{N}\frac{g_{n}(\omega)}{\pi n}e_{n}\,.
$$
Observe that in order to prove (\ref{ev1}) it suffices to establish the bound
\begin{equation}\label{ev2}
\exists\, C>0,\quad \forall\, q\geq p,\,\quad
\|S(t) f^{\omega}_{N}\|_{L^{q}(\Omega;L^p((0,2)\times \Theta))}\leq C\sqrt{q}.
\end{equation}
Indeed, if we suppose that (\ref{ev2}) holds true then by the Bienaym\'e-Tchebichev inequality, we have
$$
p(\{\omega\in \Omega\,:\, \|S(t) f^{\omega}_{N}\|_{L^p((0,2)\times \Theta)}>\lambda\})
\leq
\lambda^{-q}\|S(t) f^{\omega}_{N}\|_{L^{q}(\Omega;L^p((0,2)\times \Theta))}^{q}
\leq 
C\Big(\frac{\sqrt{q}}{\lambda}\Big)^{q}
$$
and (\ref{ev1}) follows by taking $q=\lambda^2/2$.

Let us now turn to the proof of (\ref{ev2}). Recall the general Gaussian bound (see e.g. \cite{AT,BT2})
\begin{equation}\label{ev3}
\exists\, C>0,\quad \forall\, q\geq 2,\quad\forall\, (c_n)_{n=1}^{\infty}\in l^2,\quad
\big\|\sum_{n=1}^{\infty}c_n g_{n}(\omega)\big\|_{L^q(\Omega)}
\leq C\sqrt{q}\big(\sum_{n=1}^{\infty}|c_n|^2\big)^{\frac{1}{2}}
\end{equation}
(observe that (\ref{ev3}) also follows from Lemma~4.2 of Part I).

 For $q\geq p$, using the Minkowski inequality, we can write
\begin{equation}
\|S(t)f^{\omega}_{N}\|_{L^q(\Omega;L^p((0,2)\times\Theta))}\leq 
\| S(t)f^{\omega}_{N}\|_{L^{p}((0,2)\times\Theta;L^{q}(\Omega))}\,.
\end{equation}
For fixed $t,r$, using the Gaussian bound (\ref{ev3}), we get
$$
\| S(t)f^{\omega}_{N}\|_{L^{q}(\Omega)}\leq
C\sqrt{q}\Bigl(\sum_n\Bigl|e^{-i t \pi n}\frac 1 {\pi n} e_n (r)\Bigr|^2\Bigr)^{1/2}\,. 
$$
Therefore, using that $p \geq 2$,
\begin{equation*}
\begin{aligned}
\|S(t)f^{\omega}_{N}\|_{L^q(\Omega;L^p((0,2)\times\Theta))}
&\leq 
C\sqrt{q}
\Big\|\Bigl(\sum_n\Bigl|e^{-i t \pi n}\frac {e_n(r)} n\Bigr|^2\Bigr)^{1/2}\Big\|_{L^{p} ((0,2)\times\Theta)}\\
&\leq C\sqrt{q}
\Big\|\sum_n\Bigl|\frac {e_n} n\Bigr|^2\Big\|^{1/2}_{L^{p/2} (\Theta)}\\
&\leq  C\sqrt{q}
\Bigl(\sum_n\Big\|\Bigl|\frac {e_n} n\Bigr|^2\Big\|_{L^{p/2} (\Theta)}\Bigr)^{1/2}\\
&\leq  C\sqrt{q}
\Bigl(\sum_n\Big\|\frac {e_n} n\Big\|^2_{L^{p} (\Theta)}\Bigr)^{1/2}\,.
\end{aligned}
\end{equation*}
Next (see \cite[Lemma~2.5]{AT}), we use the estimate
$$
\forall\, p\geq 2,\quad \exists\, C>0,\quad \forall\, n\geq 1,\quad \|e_n\|_{L^p(\Theta)}\leq 
Cn^{1-\frac{3}{p}}\,.
$$
This gives ($p<6$)
\begin{equation*}
\|S(t)f^{\omega}_{N}\|_{L^q(\Omega;L^p((0,2)\times\Theta))}
\leq 
C\sqrt{q}\,
\sum_n \frac 1 {n^{6/p}}\leq C \sqrt{q}
\end{equation*}
which  completes the proof of (\ref{ev2}).
This ends the proof of Proposition~\ref{largedev}.
\end{proof}
\section{Deterministic Strichartz estimates}
In this section we recall the Strichartz estimates for the free evolution 
(see~\cite[Section~4]{BT1} and \cite[Section~2]{BT2} for the proofs).
\begin{definition}
A couple of real numbers $(p,q), 2<p\leq + \infty$ is admissible if $\frac{1}{p}+\frac{1}{q}=\frac{1}{2}$.
For $T>0$, $0\leq \sigma < 1$, we define the spaces 
\begin{equation}\label{eq.espstri}
X^\sigma_T= C^0 ([-T, T]; H^\sigma(\Theta)) \cap L^p((-T, T) ; L^q(\Theta)), 
(p = \frac 2 \sigma,q) \text{ admissible} 
\end{equation}
and its dual space
\begin{equation}\label{eq.spdual}
 Y^\sigma_T= L^1 ([-T, T]; H^{-\sigma}(\Theta)) + L^{p'}((-T, T) ; L^{q'}(\Theta)), 
(p = \frac 2 \sigma,q) \text{ admissible}
\end{equation}
equipped with their natural norms ($(p',q')$ being the conjugate couple of $(p,q)$).
\end{definition} 
Remark that a simple interpolation argument gives the following statement.
\begin{lemme}\label{lem.interpo}
Assume that $0\leq \sigma <1$ and 
$$ 
\frac 1 p + \frac 3 q = \frac 3 2 -\sigma,\quad \frac{2}{\sigma}\leq p\leq  + \infty\,.
$$
Then 
\begin{equation}\label{eq.interpo}
\|f\|_{L^p([0,T]; L^q(\Theta))}\leq 
C\|f\|_{X^\sigma_T},\qquad \|f\|_{Y^\sigma_T} \leq C\|f\|_{L^{p'}([0,T]; L^{q'}(\Theta))}
\end{equation}
\end{lemme}
Recall that $S(t)=e^{-it\sqrt{- \mathbf{\Delta}}}$.
We next state several Strichartz inequalities for $S(t)$.
We refer to \cite[Section 4]{BT1} for the proof.
\begin{proposition}\label{stri}
For every $0< \sigma <1 $, every admissible couple $(p,q)$, there exists $C>0$ such that
for every $T\in]0,1]$, every $f\in H^{\frac{2}{p}}(\Theta)$ one has
\begin{equation}\label{eq.homogbis}
\|S(t)(f)\|_{X^\sigma_T}\leq 
C\|f\|_{H^{\frac{2}{p}}(\Theta)},\text{ if\, $\frac 1 p = \frac \sigma 2$}
\end{equation}
\begin{equation}\label{eq.inhomog}
\|\int_0^t S(t-\tau)\sqrt{- \mathbf{\Delta}}^{-1}(f)(\tau) d\tau\|_{X^\sigma_T}\leq C \|f\|_{Y^{1-\sigma}_T}
\end{equation} 
\begin{equation}\label{eq.inhomogbis}
\|(1- S_N)\int_0^t\sqrt{-\mathbf{\Delta}} ^{-1} 
e^{-i (t-\tau)\sqrt{ - \mathbf{\Delta}}} (f)(\tau) d\tau\|_{X^\sigma_T}
\leq C N^{\sigma-\sigma_1}\|f\|_{Y^{1- \sigma_1}_T}, \text{ if $\sigma<\sigma_1<1 $}.
\end{equation} 
\end{proposition}
\begin{remarque}
{\rm
The map $S_N$ involved in (\ref{eq.inhomogbis}) is slightly different than the corresponding one involved
in \cite{BT1}. However the proof of \cite{BT1} still works since we have that
$(1- S_N) $ is bounded from $H^{\sigma_1}$ to $H^\sigma$ with norm $\leq CN^{\sigma-\sigma_1}$.
}
\end{remarque}
We shall also make use of the next Strichartz estimate.
\begin{proposition}\label{str_pak}
Let $p$ be such that $\max(4,2\alpha)<p<6$.
Define $\sigma$ by $\sigma=\frac{3}{2}-\frac{4}{p}$.
Then there exist $C>0$ such that
for every $T\in]0,1]$, every $f\in H^{\sigma}(\Theta)$ one has
$$
\|S(t)(f)\|_{L^p([-T,T]\times \Theta)}\leq C\|f\|_{H^{\sigma}(\Theta)}\,.
$$
\end{proposition}
\begin{proof}
Let $q$ be such that $(p,q)$ is an admissible couple. Then the Sobolev
inequality and the endpoint of (\ref{eq.homogbis}) yield
\begin{eqnarray*}
\|S(t)(f)\|_{L^p([-T,T]\times \Theta)}
& \leq & C
\|(1-\mathbf{\Delta})^{\frac 1 2(\frac{3}{2}-\frac{6}{p})}S(t)(f)\|_{L^p([-T,T];L^q(\Theta))}
\\
& \leq & C\|(1-\mathbf{\Delta})^{\frac 1 2 (\frac{3}{2}-\frac{6}{p})}(f)\|_{H^{\frac{2}{p}}(\Theta)}=C\|f\|_{H^{\sigma}(\Theta)}\,.
\end{eqnarray*}
This completes the proof of Proposition~\ref{str_pak}.
\end{proof}
\section{Local well-posedness}
The problem (\ref{2}) is reduced to the integral equation
\begin{equation}\label{Duhamel}
u(t)=S(t)u_0-
i\int_{0}^{t}S(t-\tau)\sqrt{- \mathbf{\Delta}}^{-1}\big(|\Re(u(\tau))|^{\alpha}\Re(u(\tau))\big)d\tau\,.
\end{equation} 
The next statement provides bounds on the right hand-side of (\ref{Duhamel}).
\begin{proposition}\label{calculus}
For a given positive number $\alpha<3$ we choose a real number $p$ such that $\max(4,2\alpha)<p<6$.
Then we fix a real number $\sigma$ by $\sigma=\frac{3}{2}-\frac{4}{p}$.
Set $F(u)=|\Re (u)|^{\alpha}u$. Then  there exist $C>0$, $\delta >0$ such that for every $T\in]0,2]$, 
every $u_1,u_2\in X^\sigma_T$, every $v_1,v_2\in L^p((-T,T)\times\Theta)$ 
(radial with respect to the second variable) every $u_0\in H^{\sigma}(\Theta)$, 
\begin{equation}\label{i}
\big\|S(t)u_0 \big\|_{X^\sigma_T}\leq  C\|u_0\|_{H^{\sigma}(\Theta)}\, ,
\end{equation}
\begin{equation}\label{ii}
\Big\|\int_{0}^{t}S(t-\tau)\sqrt{- \mathbf{\Delta}}^{-1}
F(u_1(\tau)+v_1(\tau))d\tau\Big\|_{X^\sigma_T}\leq CT^{\delta}\big(
\|u_1\|^{\alpha+1}_{X^\sigma_T}+\|v_1\|^{\alpha+1}_{L^p_{T}L^p}\big),
\end{equation}
where $L^p_{T}L^p$ denotes the norm in $L^p((-T,T)\times\Theta)$. Moreover
\begin{multline}\label{ii'}
\Big\|(1- S_N)\int_{0}^{t}S(t-\tau)\sqrt{- \mathbf{\Delta}}^{-1}
F(u_1(\tau)+v_1(\tau))d\tau\Big\|_{X^\sigma_T}
\\
\leq CT^{\delta}N^{-\delta}
\big(\|u_1\|^{\alpha+1}_{X^\sigma_T}+\|v_1\|^{\alpha+1}_{L^p_{T}L^p}\big)\, ,
\end{multline}
\begin{multline}\label{iii}
\Big\|\int_{0}^{t}S(t-\tau)\sqrt{- \mathbf{\Delta}}^{-1}\Big(F(u_1(\tau)+v_1(\tau))-
F(u_2(\tau)+v_2(\tau))\Big)d\tau\Big\|_{X^{\sigma}_{T}}
\\
\leq CT^{\delta}\Big(\|u_1\|^{\alpha}_{X^\sigma_T}+
\|u_2\|^{\alpha}_{X^\sigma_T}+\|v_1\|_{L^p_TL^p}^{\alpha}+\|v_2\|_{L^p_TL^p}^{\alpha}
\Big)\Big(\|u_1-u_2\|_{X^\sigma_T}+\|v_1-v_2\|_{L^p_TL^p}\Big)
\end{multline}
and 
\begin{multline}\label{iv}
\Big\|\int_{0}^{t}S(t-\tau)\sqrt{- \mathbf{\Delta}}^{-1}S_{N}\Big(F(u_1(\tau)+v_1(\tau))-
F(u_2(\tau)+v_2(\tau))\Big)d\tau\Big\|_{X^{\sigma}_{T}}
\\
\leq CT^{\delta}\Big(\|u_1\|^{\alpha}_{X^\sigma_T}+
\|u_2\|^{\alpha}_{X^\sigma_T}+\|v_1\|_{L^p_TL^p}^{\alpha}+\|v_2\|_{L^p_TL^p}^{\alpha}
\Big)\Big(\|u_1-u_2\|_{X^\sigma_T}+\|v_1-v_2\|_{L^p_TL^p}\Big)
\end{multline}
\end{proposition}
\begin{proof}
Let us first observe that thanks to the assumption $p>4$, we have that
$\sigma>1/2$ and thus $p>2/\sigma$.
Estimate (\ref{i}) follows from Proposition~\ref{stri}.
Let us next show (\ref{ii'}).
Using (\ref{eq.inhomogbis}) and  Lemma~\ref{lem.interpo}
the left hand side of (\ref{ii'}) is bounded by
\begin{equation}\label{gin}
CN^{\sigma-\sigma_1}
\Big(\||\Re(u_1)|^{\alpha}\Re(u_1)\|_{L^{p'}((-T, T); L^{q'}(\Theta))}
+\||\Re(v_1)|^{\alpha}\Re(v_1)\|_{L^{p'}((-T, T); L^{q'}(\Theta))}\Big)
\end{equation}
where $\sigma_1$ (close to $\sigma$) is such that $\sigma<\sigma_1<1$ and will be fixed later
and $(p',q')$ are such that
$$
\frac{1}{p'}+\frac{3}{q'}=\frac{5}{2}+(1-\sigma_1)\,.
$$
We take $p'=q'$ and thus $\frac{1}{p'}=\frac{1}{q'}=\frac{7}{8}-\frac{\sigma_1}{4}$.
Therefore we can evaluate (\ref{gin}) by
$$
CN^{\sigma-\sigma_1}
\Big(\|u_1\|^{\alpha+1}_{L^{(\alpha+1)p'}_{T}L^{(\alpha+1)p'}}+
\|v_1\|^{\alpha+1}_{L^{(\alpha+1)p'}_{T}L^{(\alpha+1)p'}}
\Big)\,.
$$
Thanks to Lemma~\ref{lem.interpo} and the H\"older inequality, the proof of (\ref{ii'}) will be completed 
if we can provide that $(\alpha+1)p'<p$, i.e.
\begin{equation}\label{uslovie}
\frac{\alpha+1}{p}<\frac{1}{p'}=\frac{7}{8}-\frac{\sigma_1}{4}\,.
\end{equation}
Let us choose $\sigma_1$ as $\sigma_1=\sigma+\varepsilon$, where $\varepsilon>0$ is to be specified.
Thus
$$
\frac{7}{8}-\frac{\sigma_1}{4}=\frac{1}{2}+\frac{1}{p}-\frac{\varepsilon}{4}\,. 
$$
Hence (\ref{uslovie}) can be assured if we can choose $\varepsilon>0$ such that
$$
\frac{\alpha+1}{p}<\frac{1}{2}+\frac{1}{p}-\frac{\varepsilon}{4},
$$
i.e. $\frac{\varepsilon}{4}<\frac{1}{2}-\frac{\alpha}{p}$.
Thanks to the assumption $p>2\alpha$, we have that $\frac{1}{2}-\frac{\alpha}{p}>0$ 
and thus a proper choice of $\varepsilon>0$ is indeed possible.
This completes the proof of (\ref{ii'}).
The proof of \eqref{ii} is the same as the proof of (\ref{ii'}) by choosing $\sigma_1=\sigma$.
The proofs of (\ref{iii}) and (\ref{iv}) are very similar to that of (\ref{ii'}) by invoking the 
inequality
$$
\exists\, C>0,\quad\forall\, (x,y)\in\R^2,\quad
||x|^{\alpha}x-|y|^{\alpha}y|\leq C|x-y|(|x|^{\alpha}+|y|^{\alpha})\,.
$$
This completes the proof of Proposition~\ref{calculus}.
\end{proof}
As a consequence of Proposition~\ref{calculus}, we infer the following well-posedness results
for (\ref{2}).
\begin{proposition}\label{lwp}
For a given positive number $\alpha<3$ we choose a real number $p$ such that $\max(4,2\alpha)<p<6$.
Then we fix a real number $\sigma$ by $\sigma=\frac{3}{2}-\frac{4}{p}$.
There exist $C>0$, $c\in]0,1]$, $\gamma>0$ such that for every $A>0$ if we set $T=c(1+A)^{-\gamma}$ 
then for every radially symmetric $u_0$ 
satisfying $\|S(t)u_0\|_{L^p((0,2)\times\Theta)}\leq A$ there exists a unique
solution $u$ of (\ref{2}) such that
$
u(t)=S(t)u_0+v(t)
$ 
with $v\in X^{\sigma}_{T}$. Moreover $\|v\|_{X^\sigma_T}\leq CA$. In
particular, since $S(t)$ is $2$ periodic and thanks to the Strichartz
estimate of Proposition~\ref{str_pak}, 
$$
\sup_{t\in[-T,T]}\|S(\tau)u(t)\|_{L^p(\tau\in (0,2);L^p(\Theta))}\leq CA\,.
$$
In addition, if $u_0\in H^s(\Theta)$ (and thus $s<\sigma$) then
$$
\|u(t)\|_{H^s(\Theta)}\leq
\|S(t)u_0\|_{H^s(\Theta)}+\|v(t)\|_{H^s(\Theta)}\leq
\|u_0\|_{H^s(\Theta)}+CA\,.
$$
\end{proposition}
\begin{remarque}
{\rm
Thanks to Proposition~\ref{largedev} the data
$$
f^{\omega}=\sum_{n=1}^{\infty}\frac{g_{n}(\omega)}{\pi n}e_{n}
$$
satisfies the assumption  $\|S(t)f^{\omega}\|_{L^p((0,2)\times\Theta)}<\infty$,
almost surely in $\omega$.
Therefore, despite the fact that $f^{\omega}$ is essentially in $H^{1/2}$ and
not more regular, and thus supercritical for (\ref{2}) for $2<\alpha<3$,
Proposition~\ref{lwp} establishes a local well-posedness theory for data
$f^{\omega}$ almost surely in $\omega$. We refer to Part I (cf. \cite{BT2})
for a general local well-posedness theory for the cubic wave equation posed on
a compact manifold with random initial data.
}
\end{remarque}
\begin{proof}[Proof of Proposition~\ref{lwp}.]
If we write $u(t)=S(t)u_0+v(t)$ 
then $v(0)=0$ and $v$ solves the problem
$$
(i\partial_{t}-\sqrt{ - \mathbf{\Delta}})v-\sqrt{ - \mathbf{\Delta}}^{-1}
\big(|\Re(S(t)u_0+v)|^{\alpha}\Re(S(t)u_0+v)\big)=0
$$
with corresponding integral equation
$$
v(t) = - i\int_{0}^{t}
S(t-\tau)\sqrt{-\mathbf{\Delta}}^{-1}\big(|\Re(S(\tau)u_0+v(\tau))|^{\alpha}\Re(S(\tau)u_0+v(\tau))\big)
d\tau
\equiv  K_{u_0}(v)\,.
$$
Using (\ref{ii}) and (\ref{iii}) of Proposition~\ref{calculus}, we infer that for $u_0$ such that for $T\in ]0,2]$,
$\|S(t)u_0\|_{L^p((0,2)\times\Theta)}\leq A$,
$$
\| K_{u_0}(v)\|_{X^{\sigma}_{T}}\leq
CT^{\delta}A^{\alpha+1}+CT^{\delta}\|v\|_{X^{\sigma}_{T}}^{\alpha+1}
$$
and
$$
\| K_{u_0}(v_1)-K_{u_0}(v_2)\|_{X^{\sigma}_{T}}\leq
CT^{\delta}\|v_1-v_2\|_{X^{\sigma}_{T}}
(A^{\alpha}+\|v_1\|_{X^{\sigma}_{T}}^{\alpha}+\|v_2\|_{X^{\sigma}_{T}}^{\alpha})\,.
$$
Proposition~\ref{lwp} follows by applying the contraction mapping principle
to the nonlinear map $K_{u_0}$ on the ball of radius $A$ of $X_{T}^{\sigma}$ (centered
at the origin) with $T=c(1+A)^{-\gamma}$ for a suitable choice of $c\ll 1$ and
$\gamma\gg 1$. 
\end{proof}
Thanks to (\ref{iv}) and the fact that $S_N$ is (uniformly with respect to $N$) bounded on $X^{\sigma}_{T}$
(see Lemma~\ref{psdo}), we can apply the argument of the proof of Proposition~\ref{lwp}
to obtain a well-posedness
in the context of (\ref{pak_N}) with bounds independent of $N$. 
\begin{proposition}\label{lwpbis}
For a given positive number $\alpha<3$ we choose a real number $p$ such that $\max(4,2\alpha)<p<6$.
Then we fix a real number $\sigma$ by $\sigma=\frac{3}{2}-\frac{4}{p}$. For
every $N\geq 1$
there exist $C>0$, $c\in]0,1]$, $\gamma>0$ such that for every 
$A>0$ if we set $T=c(1+A)^{-\gamma}$ then for every $u_0\in E_{N}$ 
satisfying 
$\|S(t)u_0\|_{L^p((0,2)\times\Theta)}\leq A$
there exists a unique solution $u$ of 
(\ref{pak_N}) 
such that
$
u(t)=S(t)u_0+v(t)
$ 
with $v\in X^{\sigma}_{T}$. Moreover $\|v\|_{X^\sigma_T}\leq CA$ and in
particular, 
$$
\sup_{t\in[-T,T]}\|S(\tau)u(t)\|_{L^p(\tau\in (0,2);L^p(\Theta))}\leq CA\,.
$$
In addition,
$$
\|u(t)\|_{H^s(\Theta)}\leq
\|S(t)u_0\|_{H^s(\Theta)}+\|v(t)\|_{H^s(\Theta)}\leq
\|u_0\|_{H^s(\Theta)}+CA\,.
$$

\end{proposition}
\section{Global existence for (\ref{2}) on a set of full $\rho$ measure}
Let us denote by $\Phi_{N}(t):E_{N}\rightarrow E_{N}$, $t\in\R$ the flow of (\ref{pak_N}) defined in
Proposition~\ref{liouville}.
In the next proposition, we obtain a crucial long time bound for the solutions
of (\ref{pak_N}) (a similar argument was already performed in \cite{Tz1,Tz2,BT1}).
\begin{proposition}\label{longtime}
Let us fix $p$ such that $2\alpha<p<6$. 
Then for every integer $i\geq 1$, every integer $N\geq 1$, there exists 
a $\rho_{N}$ measurable set $\Sigma_{N}^{i}\subset E_{N}$ such that
$
\rho_{N}(E_{N}\backslash \Sigma_{N}^{i})\leq 2^{-i}
$
and there exists a constant $C$ such that for every $i, N\in\N$, every $u_0\in \Sigma_{N}^{i}$, every $t\in\R$, 
\begin{equation}\label{uniform}
\|S(\tau)(\Phi_{N}(t)(u_0))\|_{L^p(\tau\in (0,2);L^p(\Theta))}
+
\|\Phi_{N}(t)(u_0)\|_{H^{s}(\Theta)}\leq C(i+\log(1+|t|))^{\frac{1}{2}}\,.
\end{equation}
\end{proposition}
\begin{proof}
For $i,j$ integers $\geq 1$, we set
$$
B_{N}^{i,j}(D)\equiv
\big\{u\in E_{N}\,:\,
\|S(t)u\|_{L^p((0,2)\times\Theta)}
+\|u\|_{H^{s}(\Theta)}\leq D(i+j)^{\frac{1}{2}}\big\},
$$
where the number $D\gg 1$ (independent of $i,j,N$) will be fixed later. 
Thanks to Lemma~\ref{lem.gauss}, Proposition~\ref{largedev}, we have that
\begin{equation}\label{residuel}
\rho_{N}(E_{N}\backslash B_{N}^{i,j}(D))\leq Ce^{-cD^2(i+j)}\,.
\end{equation}
Thanks to Proposition~\ref{lwpbis}, there exist $c>0$, $C>0$, $\gamma>0$ only depending on $\alpha$  
such that if we set $\tau \equiv c D^{-\gamma}(i+j)^{-\gamma/2}$ then for every $t\in[-\tau,\tau]$,
\begin{equation}\label{preser}
\Phi_{N}(t)\big(B_{N}^{i,j}(D)\big)\subset 
\Big\{u\in E_{N}\,:\,
\|S(t)u\|_{L^p((0,2)\times\Theta)}+\|u\|_{H^{s}(\Theta)}\leq C\,D(i+j)^{\frac{1}{2}}\}\, .\end{equation}
Next, we set
$$
\Sigma_{N}^{i,j}(D)\equiv
\bigcap_{k=-[2^{j}/\tau]}^{[2^{j}/\tau]}\Phi_{N}(-k\tau)(B_{N}^{i,j}(D))\, ,
$$
where $[2^{j}/\tau]$ stays for the integer part of $2^{j}/\tau$. 
Using the invariance of the measure $\rho_{N}$ by the flow $\Phi_{N}$ (Proposition~\ref{liouville}), 
we can write
\begin{equation*}
\rho_{N}(E_{N}\backslash\Sigma_{N}^{i,j}(D))
\leq 
(2[2^{j}/\tau]+1)\rho_{N}(E_N\backslash B_{N}^{i,j}(D))
\leq 
C2^{j}D^{\gamma}(i+j)^{\gamma/2}\rho_{N}(E_N\backslash B_{N}^{i,j}(D))\,.
\end{equation*}
Using (\ref{residuel}), we now deduce 
\begin{equation}\label{zvez}
\rho_{N}(E_{N}\backslash\Sigma_{N}^{i,j}(D))\leq
C2^{j}D^{\gamma}(i+j)^{\gamma/2}e^{-cD^2(i+j)}\leq 2^{-(i+j)},
\end{equation}
provided $D\gg 1$, independently of $i,j,N$.
Thanks to (\ref{preser}), we obtain that for
$u_0\in\Sigma_{N}^{i,j}(D)$, the solution of (\ref{pak_N}) with data $u_0$ satisfies
\begin{equation}\label{jjj1}
\big\|S(\tau)\big(\Phi_{N}(t)(u_0)\big)\big\|_{L^p(\tau\in (0,2);L^p(\Theta))}
+
\|\Phi_{N}(t)(u_0)\|_{H^{s}(\Theta)}
\leq CD(i+j)^{\frac{1}{2}},\quad |t|\leq 2^{j}\,.
\end{equation}
Indeed, for $|t|\leq 2^{j}$, we may find an integer $k\in [-[2^{j}/\tau],[2^{j}/\tau]]$ and 
$\tau_1\in [-\tau,\tau]$ so that $t=k\tau+\tau_1$ and thus 
$u(t)=\Phi_{N}(\tau_1)\big(\Phi_{N}(k\tau)(u_0)\big)$.
Since $u_0\in\Sigma_{N}^{i,j}(D)$ implies that $\Phi_{N}(k\tau)(u_0)\in 
B_{N}^{i,j}(D)$, we may apply (\ref{preser}) and arrive at (\ref{jjj1}).
Next, we set
$$
\Sigma_{N}^{i}=\bigcap_{j= 1}^{\infty}\Sigma_{N}^{i,j}(D)\,.
$$
Thanks to (\ref{zvez}),
$
\rho_{N}(E_{N}\backslash \Sigma_{N}^{i})\leq 2^{-i}\,.
$
In addition, using (\ref{jjj1}), we get that there exists $C$ such that for every $i$, every
$N$, every $u_0\in \Sigma_{N}^{i}$, every $t\in \R$,
\begin{equation}\label{estmi}
\big\|S(\tau)\big(\Phi_{N}(t)(u_0)\big)\big\|_{L^p(\tau\in (0,2);L^p(\Theta))}+\|\Phi_{N}(t)(u_0)\|_{H^{s}(\Theta)}
\leq C(i+\log(1+|t|))^{\frac{1}{2}}\,.
\end{equation}
Indeed for $t\in \R$ there exists $j\in\N$ such that $2^{j-1}\leq 1+|t|\leq 2^j$ and we apply 
(\ref{jjj1}) with this $j$.
This completes the proof of Proposition~\ref{longtime}.
\end{proof}
For  integers $i\geq 1$ and $N \geq 1$, we define the
cylindrical sets
$$
\tilde{\Sigma}_{N}^{i}\equiv\big\{u\in H^{s}(\Theta)\,:\, S_{N}(u)\in \Sigma_{N}^{i}\big\}.
$$
Next, for an integer $i\geq 1$, we set
\begin{equation*}
\Sigma^{i}\equiv \limsup_{N\rightarrow\infty}\tilde{\Sigma}_{N}^{i}
\equiv\bigcap_{N=
  1}^{\infty}\bigcup_{N_1=N}^{\infty}\tilde{\Sigma}_{N_1}^{i}\,.
\end{equation*}
Using Fatou's lemma, we get
\begin{equation}\label{kr2}
\rho(\limsup_{N\rightarrow\infty}\tilde{\Sigma}_{N}^{i})
 \geq 
\limsup_{N\rightarrow \infty}\rho(\tilde{\Sigma}_{N}^{i})\,.
\end{equation}
We have that
$$
\rho(\tilde{\Sigma}_{N}^{i})=\int_{\tilde{\Sigma}_{N}^{i}}f(u)d\mu(u)
$$
and
$$
\rho_{N}(\Sigma_{N}^{i})=\int_{\Sigma_{N}^{i}}f_{N}(u)d\mu_{N}(u)=\int_{\tilde{\Sigma}_{N}^{i}}f_{N}(u)d\mu(u)
$$
where $f$ and $f_N$ are defined by (\ref{f_N}).
Therefore, thanks to Lemma~\ref{deg}, we get
$$
\lim_{N\rightarrow \infty}\big((\rho(\tilde{\Sigma}_{N}^{i})-\rho_{N}(\Sigma_{N}^{i})\big)=0\,.
$$
Therefore, using  Proposition~\ref{longtime} and (\ref{q}), we obtain
\begin{equation}\label{kr4}
\limsup_{N\rightarrow \infty}\rho(\tilde{\Sigma}_{N}^{i})
=
\limsup_{N\rightarrow \infty}\rho_{N}(\Sigma_{N}^{i})
\geq
\limsup_{N\rightarrow \infty}\big(\rho_{N}(E_{N})-2^{-i}\big)
=
\rho\big(H^{s}(\Theta)\big)-2^{-i}.
\end{equation}
Collecting (\ref{kr2}) and (\ref{kr4}), we arrive at
$$
\rho(\Sigma^{i})  \geq \rho\big(H^{s}(\Theta)\big)-2^{-i}.
$$
Now, we set
$$
\Sigma\equiv\bigcup_{i\geq 1}\Sigma^{i}\,.
$$
Thus $\Sigma$ is of full $\rho$ measure.
It turns out that one has global existence for $u_0\in \Sigma$. 
\begin{proposition}\label{global_existence}
Choose a real number $p$ such that $\max(4,2\alpha)<p<6$
and then a real number $\sigma$ by $\sigma=\frac{3}{2}-\frac{4}{p}$
(so that we are in the scope of the applicability of
Propositions~\ref{lwp}, \ref{lwpbis}, \ref{longtime}).
Let us fix $i\in\N$. Then for every $u_0\in \Sigma^{i}$, the local solution $u$ 
of (\ref{2}) given by Proposition~\ref{lwp} is globally defined.
In addition there exists $C>0$ such that for every $u_0\in \Sigma^{i}$,
\begin{equation}\label{growth}
\forall t \in \mathbb{R}, \quad\|u(t)\|_{H^s(\Theta)}
+
\|S(\tau)(u(t))\|_{L^p(\tau\in (0,2);L^p(\Theta))}
\leq C(i+\log(1+|t|))^{\frac{1}{2}}\,.
\end{equation}
\end{proposition}
\begin{proof}
Let $u_0\in \Sigma^{i}$.
Then there exists $N_k\rightarrow\infty$ such that $S_{N_k}(u_0)\in
\Sigma_{N_k}^{i}$. Set $u_{0,k}\equiv S_{N_k}(u_0)$. Then $u_{0,k}$
is a sequence such that
$$
\lim_{k\rightarrow\infty}
\|u_0-u_{0,k}\|_{H^s(\Theta)}=0.
$$
Furthermore, thanks to~(\ref{uniform}), we have
$$ \|S(t)u_{0,k}\|_{L^p((0,2)\times\Theta)}\leq C i.$$
After possibly extracting a subsequence, we have that $S(t)u_{0,k}$ converges in $L^p$ for the weak topology to a function $g\in L^p((0,2)\times\Theta)$. But, as $S(t)u_{0,k}$ converges in $\mathcal{D}'$ to $S(t)u_{0}$, we deduce that $S(t)u_{0}=g \in  L^p((0,2)\times\Theta)$. But, thanks to
Lemma~\ref{psdo}, the family $(S_{N_k})_{k=1}^\infty$ is uniformly bounded on $  L^p((0,2)\times\Theta)$, and
$$ 
\forall\, g \in L^p((0,2)\times\Theta),\quad 
\lim_{k\rightarrow + \infty}\|g- S_{N_k}g\|_{L^p((0,2)\times\Theta)} =0.
$$
Indeed, it is true if $g \in C^\infty_0( (0,2)\times\Theta)$ and follows for general $g$ by density. As a consequence, we deduce
\begin{equation}\label{car}
\lim_{k\rightarrow\infty}
 \|S(t)(u_0-u_{0,k})\|_{L^p((0,2)\times\Theta)}+\|u_0-u_{0,k}\|_{H^s(\Theta)}=0.
\end{equation}

Let us fix $T>0$. Our aim is to extend the solution 
of (\ref{2}) given by Proposition~\ref{lwp} to the interval $[-T,T]$.
Using Proposition~\ref{longtime}, we have that there exists a constant $C$ such that for every $k\in\N$,
every $t\in\R$,
\begin{equation}\label{ant}
\|S(\tau)(\Phi_{N_k}(t)(u_{0,k}))\|_{L^p(\tau\in (0,2);L^p(\Theta))}
+\|\Phi_{N_k}(t)(u_{0,k})\|_{H^{s}(\Theta)}\leq C(i+\log(1+|t|))^{\frac{1}{2}}\,.
\end{equation}
To prove~\eqref{growth}, we are going to pass to the limit $k \rightarrow + \infty$ in ~\eqref{ant}. If we set $u_{N_k}(t)\equiv  \Phi_{N_k}(t)(u_{0,k})$ and 
$\Lambda\equiv C(i+\log(1+T))^{\frac{1}{2}}$, we have the bound
\begin{equation}\label{david1}
\|S(\tau)(u_{N_k}(t))\|_{L^p(\tau\in (0,2);L^p(\Theta))}
+\|u_{N_k}(t)\|_{H^{s}(\Theta)}\leq \Lambda,\quad \forall\,|t|\leq T,\,\forall\, k\in\N.
\end{equation}
In particular 
\begin{equation}\label{u0}
\|S(\tau)(u_{0})\|_{L^p(\tau\in (0,2);L^p(\Theta))}
+\|u_0\|_{H^{s}(\Theta)}\leq\Lambda
\end{equation}
(apply (\ref{david1}) with $t=0$ and let $k\rightarrow\infty$).
Let $\tau>0$ be the local existence time for (\ref{2}), provided by
Proposition~\ref{lwp} for $A=\Lambda$. 
Recall that we can assume $\tau=c(1+\Lambda)^{-\gamma}$
for some $c>0$, $\gamma>0$ depending only on $p$. 
We can also assume that $T>\tau$.
Denote by $u(t)$ the solution of (\ref{2}) with data $u_0$ on the time
interval $[-\tau,\tau]$. 
Define $v$ by $u(t)=S(t)(u_0)+v(t)$.
Thanks to (\ref{u0}) and Proposition~\ref{lwp}, we have that
\begin{equation}\label{ant1}
\|v\|_{X^{\sigma}_{\tau}}+\|u(t)\|_{H^s(\Theta)}+\|S(\tau_1)(u(t))\|_{L^p(\tau_1\in(0,2);L^p(\Theta))}\leq
C\Lambda,\quad t\in [-\tau,\tau],
\end{equation}
where $C$ depends only on $p$.
Next we define
$v_{N_k}(t)$ by $u_{N_k}(t)= S(t)(u_{0,k})+v_{N_k}(t)$.
Thanks to (\ref{david1}) and Proposition~\ref{lwpbis}, we have that
\begin{equation}\label{ant2}
\|v_{N_k}\|_{X^{\sigma}_{\tau}}+
\|u_{N_k}(t)\|_{H^s(\Theta)}+\|S(\tau_1)(u_{N_k}(t))\|_{L^p(\tau_1\in(0,2);L^p(\Theta))}\leq
C\Lambda,\quad t\in [-\tau,\tau]\,.
\end{equation}
We have that
$
w_{N_k}\equiv v-v_{N_k}
$
solves the equation
\begin{equation}\label{eqnv}
(i\partial_{t}-\sqrt{- \mathbf{\Delta}}) w_{N_k} =
\sqrt{- \mathbf{\Delta}}^{-1}\Big(F(u)-S_{N_k}(F(S_{N_k}(u_{N_k})))\Big),
\quad  w_{N_k}|_{t=0}=0,
\end{equation}
where $F(u)=|\Re(u)|^{\alpha}\Re(u)$.
Next, we write
$$
F(u)-S_{N_k}(F(S_{N_k}(u_{N_k})))=S_{N_k}\big(F(u)-F(S_{N_k}(u_{N_k}))\big)+(1-S_{N_k})F(u).
$$
Therefore
\begin{multline}\label{back}
w_{N_k}(t)=
-i\int_{0}^{t}S(t-\tau)\sqrt{- \mathbf{\Delta}}^{-1}S_{N_k}\big(F(u(\tau))-F(S_{N_k}(u_{N_k})(\tau))\big)d\tau
\\
-i\int_{0}^{t}S(t-\tau)\sqrt{- \mathbf{\Delta}}^{-1}(1-S_{N_k})F(u(\tau))d\tau.
\end{multline}
Using Proposition~\ref{calculus}, we obtain that there exist $C>0$ and
$\theta, \delta>0$ depending only on $p$ such that one has the bound
\begin{equation*}
\|(1- S_{N_k})\int_{0}^{t}S(t-\tau)\sqrt{-\mathbf{\Delta}}^{-1}
F(u(\tau))d\tau\|_{X^{\sigma}_{\tau}}\leq C \tau^{\theta}N_k^{- \delta}
\big(\|S(t)(u_0)\|_{L^p((-\tau,\tau)\times\Theta)}+\|v\|_{X^{\sigma}_{\tau}}^{\alpha}\big). 
\end{equation*}
Another use of Proposition~\ref{calculus} yields
\begin{multline*}
\|\int_{0}^{t}S(t-\tau)\sqrt{- \mathbf{\Delta}}^{-1}S_{N_k}\big(F(u(\tau))-
F(S_{N_k}(u_{N_k})(\tau))\big)d\tau\|_{X^{\sigma}_{\tau}}
\\
\leq C\tau^{\theta}
\big(\|S(t)(u_0-S_{N_k}(u_{0,k}))\|_{L^p((-\tau,\tau)\times\Theta)}+\|v-S_{N_k}(v_{N_k})\|_{X^{\sigma}_{\tau}}\big)
\times 
\\
\big(
\|S(t)(u_0)\|_{L^p((-\tau,\tau)\times\Theta)}^{\alpha}+\|S(t)(u_{0,k})\|_{L^p((-\tau,\tau)\times\Theta)}^{\alpha}
+\|v_{N_k}\|_{X^{\sigma}_{\tau}}^{\alpha}+\|v\|_{X^{\sigma}_{\tau}}^{\alpha}\big).
\end{multline*}
Collecting the last two bounds (\ref{david1}), (\ref{u0}), (\ref{ant1}), (\ref{ant2}),
coming back to (\ref{back}) yields
\begin{equation*}
\|w_{N_k}\|_{X^{\sigma}_{\tau}} \leq
C\tau^{\theta}(1+\Lambda)^{\alpha}\|w_{N_k}\|_{X^{\sigma}_{\tau}} + o(1)_{k\rightarrow + \infty}\,.
\end{equation*}
Recall that $\tau=c(1+\Lambda)^{-\gamma}$, where $c>0$ and $\gamma>0$ are depending only on 
$p$. In the last estimate the constants $C$ and $\theta$ also depend only
on $p$. 
Therefore, if we assume that $\gamma>\alpha/\theta$ then the restriction on $\gamma$
remains to depend only on $p$. Similarly, if we assume that $c$ is so small that
$
C\tau^{\theta}(1+\Lambda)^{\alpha}\leq
Cc^{\theta}(1+\Lambda)^{-\gamma\theta}(1+\Lambda)^{\alpha}
\leq Cc^{\theta}
<1/2
$
then the smallness restriction on $c$ remains to depend only on $p$. 
Therefore, we have that after possibly slightly modifying the values of $c$ and $\gamma$ 
(keeping $c$ and $\gamma$ independent of $N_k$)
in the definition of $\tau$ that
\begin{eqnarray}\label{eq.est}
\lim_{k\rightarrow\infty}\|w_{N_k}\|_{X^{\sigma}_{\tau}}=0,
\end{eqnarray}
where $ \tau= c(1+\Lambda)^{-\gamma}$
and the constants $c$ and $\gamma$ depend only on $p$.
Hence
\begin{equation}\label{limit2}
\lim_{k\rightarrow\infty}\|u-u_{N_k}-S(t)(u_0-u_{0,k})\|_{X^{\sigma}_{\tau}}=0\,.
\end{equation}
Coming back to (\ref{car}), we obtain that
\begin{equation}\label{iter1}
\lim_{k\rightarrow\infty}\|u(\tau)-u_{N_k}(\tau)\|_{H^s(\Theta)}=0\,.
\end{equation}
Moreover combining (\ref{limit2}) with (\ref{car}) and the Strichartz
inequality of Proposition~\ref{str_pak} yields
\begin{equation}\label{iter2}
\lim_{k\rightarrow\infty}
\|S(\tau_1)(u(\tau)-u_{N_k}(\tau))\|_{L^p(\tau_1\in(0,2);L^p(\Theta))}=0.
\end{equation}
As a consequence of (\ref{iter1}), (\ref{iter2}) and (\ref{david1}), we infer that
\begin{equation}\label{iter3}
\|u(\tau)\|_{H^s(\Theta)}+\|S(\tau_1)(u(\tau)\|_{L^p(\tau_1\in(0,2);L^p(\Theta))}\leq
\Lambda\,.
\end{equation}
Thanks to (\ref{iter1}), (\ref{iter2}) and (\ref{iter3}) we can repeat the argument on 
$(\tau,2\tau)$, $(2\tau, 3\tau)$, ...$([\frac T \tau] \tau,([\frac T \tau]+1) \tau)$
(and similarly for negative times), giving existence up to the time $T$ (which
was an arbitrary number) and (\ref{growth}).
This completes the proof of Proposition~\ref{global_existence}.
\end{proof}
Therefore we solved globally the problem
(\ref{2}) on a set of full $\rho$ measure.
This completes the proof of Theorem~\ref{thm1}.
\begin{remarque}
It is likely that as in \cite{BT1}, where the easier sub-critical problem is studied, we may further push
the analysis in order to prove that the measure $\rho$ is indeed
invariant under the flow of (\ref{2}) established by  Theorem~\ref{thm1}.
We decided not to pursue this issue here since our main concern in the present
paper is to establish random data Cauchy theory for supercritical problems.
We refer to \cite{Bo1,Bo2,KS,Tz1,Tz2,Zh} for results concerning the existence of
invariant Gibbs measures in the closely related context of the Nonlinear
Schr\"odinger equation.
\end{remarque}


\begin{thebibliography}{10}
\bibitem{AT} A. Ayache, N. Tzvetkov, {\it $L^p$ properties of Gaussian random series},
to appear in Trans. AMS.
%
\bibitem{Bo1}
J.~Bourgain, {\it Periodic nonlinear Schr\"odinger equation and invariant measures}, 
Comm. Math. Phys. 166 (1994) 1-26.
%
\bibitem{Bo2}
J.~Bourgain, {\it Invariant measures for the 2D-defocusing nonlinear Schr\"odinger equation}, 
Comm. Math. Phys. 176 (1996) 421-445.
%
\bibitem{BGT}
N.~Burq, P.~G{\'e}rard, N.~Tzvetkov,
{\it {S}trichartz inequalities and the non linear {S}chr\"odinger equation on compact manifolds}, 
Amer. J. of Math., 126 (2004) 569-605.
%
\bibitem{BT1} N.~Burq, N.~Tzvetkov, {\it Invariant measure for the three dimensional nonlinear wave equation}, 
Preprint~2007, {\em http://arxiv.org/abs/0707.1445}.
%
\bibitem{BT2} N.~Burq, N.~Tzvetkov, {\it 
Random data Cauchy theory for supercritical wave equations I : local existence theory}, Preprint~2007, {\em http://arxiv.org/abs/0707.1447}. 
%
\bibitem{KS} S.~Kuksin, A.~Shirikyan,
{\it Randomly forced CGL equation : stationary measures and the inviscid limit}, J. Phys A 37 (2004) 1-18.
%
\bibitem{LionsMagenes} J.L.~Lions, Magenes {\it Probl\`emes aux limites non homog\`enes et applications}, 
Dunod, Paris 1968.
%
\bibitem{Sz} J.~Szeftel, {\it
Propagation et r\'eflexion des singularit\'es pour l'\'equation de Schr\"odinger non lin\'eaire
}, Ann. Inst. Fourier, 55 (2005) 573--671. 
%
\bibitem{Tz1} N.~Tzvetkov, {\it Invariant measures for the Nonlinear Schr\"odinger equation on the disc},
Dynamics of PDE 3 (2006) 111-160.
%
\bibitem{Tz2} N.~Tzvetkov, {\it Invariant measures for the defocusing NLS}, Preprint~2007.
%
\bibitem{Zh} P. Zhidkov, {\it KdV and nonlinear Schr\"odinger equations :
Qualitative theory}, Lecture Notes in Mathematics 1756, Springer 2001.
\end{thebibliography}
\end{document}